\documentclass[12pt]{article}
\usepackage{latexsym}
\usepackage{amsfonts}
\usepackage{amsmath}
\usepackage{amssymb}
\usepackage{graphicx}
\usepackage{latexsym}
\usepackage{amsfonts}
\usepackage{graphicx}
\usepackage{psfrag}

\input{tcilatex}
\begin{document}

\qquad 

\qquad 

\thispagestyle{empty}

\begin{center}
{\Large \textbf{Eigenvalues and Diagonal Elements}}

\vskip0.2in$^{1}$Rajendra Bhatia and $^{2}$Rajesh Sharma

$^{1}$Ashoka University, Sonepat, Haryana, 131029, India

$^{2}$Department of Mathematics \& Statistics, H.P. University Shimla-5,
India
\end{center}

\vskip0.5in \noindent \textbf{Abstract. }A basic theorem in linear algebra
says that if the eigenvalues and the diagonal entries of a Hermitian matrix
are ordered as $\lambda _{1}\leq \lambda _{2}\leq \cdots \leq \lambda _{n}$
and $a_{1}\leq a_{2}\leq \cdots \leq a_{n},$ respectively, then $\lambda
_{1}\leq a_{1}.$ We show that for some special classes of Hermitian matrices
this can be extended to inequalities of the form $\lambda _{k}\leq a_{2k-1},$
$k=1,2,...,\lceil \frac{n}{2}\rceil .$


\bigskip

\pagebreak

\setcounter{equation}{0}Let $A$ be an $n\times n$ complex Hermitian matrix.
The eigenvalues and the diagonal entries of $A$ are real numbers, and we
enumerate them in increasing order as%
\begin{equation*}
\lambda _{1}\leq \lambda _{2}\leq \cdots \leq \lambda _{n},
\end{equation*}%
and 
\begin{equation*}
a_{1}\leq a_{2}\leq \cdots \leq a_{n},
\end{equation*}%
respectively. Various inequalities relating these two $n$-tuples are known
and are much used in matrix analysis. For example, we have%
\begin{equation}
\lambda _{1}\leq a_{1}\text{ \ \ and\ \ }\lambda _{n}\geq a_{n}.  \tag{1}
\end{equation}%
These are subsumed in the majorization relations due to I. Schur: for $1\leq
k\leq n$

\begin{equation}
\underset{j=1}{\overset{k}{\sum }}\lambda _{j}\leq \underset{j=1}{\overset{k}%
{\sum }}a_{j},  \tag{2}
\end{equation}%
with equality when $k=n$. This is a complete characterization of two $n$%
-tuples that could be the eigenvalues and diagonal entries of a Hermitian
matrix. In general, there are no further relations between individual $%
\lambda _{j}$ and $a_{k}.$ However, for large and interesting subsets of
Hermitian matrices, it might be possible to find such extra relations. In 
\cite{1} the authors consider eigenvalues of matrices associated with
graphs. Let $G$ be a simple weighted graph on $n$ vertices and let $A$ be
the signless Laplacian matrix associated with $G.$ Then, it is shown in \cite%
{1} that $\lambda _{2}\leq a_{3}.$ This result is extended to other classes
in \cite{3}. One of these is the class $\mathcal{P}$ of Hermitian matrices
whose off-diagonal entries are nonnegative. (In particular, this includes
symmetric entrywise nonnegative matrices.) It is shown in \cite{3} that if $%
A\in \mathcal{P},$ then $\lambda _{2}\leq a_{3}.$

In this note we consider, in addition the class $\mathcal{P}$, another class 
$\mathcal{I}$ consisting of Hermitian matrices all whose off-diagonal
entries are purely imaginary. We show that the inequality $\lambda _{2}\leq
a_{3}$ is valid for $A\in \mathcal{I}$ as well. The proof we give works for
both the classes $\mathcal{P}$ and $\mathcal{I}$. Then we show that much
more is true for the class $\mathcal{I}.$ We show that in this case the
inequality $\lambda _{n-1}\geq a_{n-2}$ also holds. Further, for all $1\leq
k\leq \lceil \frac{n}{2}\rceil $ we have $\lambda _{k}\leq a_{2k-1}.$ We
construct examples to show that neither of these results is true for the
class $\mathcal{P}.$\vskip0.2in\noindent \textbf{Theorem 1}. Let $A$ be an $%
n\times n$ Hermitian matrix whose off-diagonal entries are either all
nonnegative real numbers or all purely imaginary numbers. Then

\begin{equation}
\lambda _{2}\leq a_{3}.  \tag{3}
\end{equation}%
In case the off-diagonal entries are all purely imaginary, we also have

\begin{equation}
\lambda _{n-1}\geq a_{n-2}.  \tag{4}
\end{equation}

For the second class of matrices in Theorem 1, we can go further:\vskip%
0.2in\noindent \textbf{Theorem 2}. Let $A$ be an $n\times n$ Hermitian
matrix whose off-diagonal entries are all purely imaginary. Then, for $1\leq
k\leq \lceil \frac{n}{2}\rceil ,$%
\begin{equation}
\lambda _{k}\leq a_{2k-1}\text{ \ \ and \ \ }\lambda _{n-k+1}\geq a_{n-2k+2}.
\tag{5}
\end{equation}

We remark that in both (1) and (5) the second inequality follows from the
first by considering $-A$ in place of $A$. Similarly (4) follows from (3).
The argument cannot be used for the class $\mathcal{P}$.

Our proofs rely upon two basic theorems of matrix analysis. Let $\lambda
_{j}\left( A\right) $, $1\leq j\leq n$, denote the eigenvalues of a
Hermitian matrix enumerated in the increasing order. Weyl's inequality says
that if $A$ and $B$ are two $n\times n$ Hermitian matrices, then%
\begin{equation}
\lambda _{j}\left( A+B\right) \leq \lambda _{j}\left( A\right) +\lambda
_{n}\left( B\right) ,\text{ \ \ \ }1\leq j\leq n.  \tag{6}
\end{equation}%
Cauchy's interlacing principle says that if $A_{r}$ is an $r\times r$
principal submatrix of $A$, then

\begin{equation}
\lambda _{j}\left( A\right) \leq \lambda _{j}\left( A_{r}\right) ,\text{ \ \
\ }1\leq j\leq r.  \tag{7}
\end{equation}%
See Chapter III of \cite{2} for this and other facts used here.\vskip%
0.2in\noindent \textbf{Proof of Theorem 1.} If $P$ is a permutation matrix,
then the increasingly ordered eigenvalues and diagonal entries of $PAP^{T}$
are the same as those of $A$. So, for simplicity, we may assume that the
diagonal entries of $A$ are in increasing order. Let

\begin{equation*}
A_{3}=\left[ 
\begin{array}{ccc}
a_{11} & a_{12} & a_{12} \\ 
\overline{a_{12}} & a_{22} & a_{23} \\ 
\overline{a_{13}} & \overline{a_{23}} & a_{33}%
\end{array}%
\right]
\end{equation*}%
be the top-left $3\times 3$ submatrix of $A$. (Note $a_{jj}=a_{j}$ is our
notation.) Decompose

\begin{equation}
A_{3}=D_{3}+M_{3}  \tag{8}
\end{equation}%
where $D_{3}$ is the diagonal part and $M_{3}$ the off-diagonal part of $%
A_{3}$. By Weyl's inequality%
\begin{equation}
\lambda _{2}\left( A_{3}\right) \leq \lambda _{2}\left( M_{3}\right)
+\lambda _{3}\left( D_{3}\right) =\lambda _{2}\left( M_{3}\right) +a_{3}. 
\tag{9}
\end{equation}%
Note that $\det M_{3}=2\func{Re}a_{12}a_{23}\overline{a_{13}}$. So, under
the hypothesis of Theorem 1, $\det M_{3}\geq 0.$ We also have tr$M_{3}=0$.
These two conditions imply that we must have $\lambda _{2}\left(
M_{3}\right) \leq 0$. For, if $\lambda _{3}\left( M_{3}\right) \geq \lambda
_{2}\left( M_{3}\right) >0$, then the condition tr$M_{3}=0$ forces $\lambda
_{1}\left( M_{3}\right) $ to be negative. But this is impossible if $\det
M_{3}\geq 0.$ So, from (9) we see that $\lambda _{2}\left( A_{3}\right) \leq
a_{3}.$ Then, by the interlacing principle (7), we have $\lambda _{2}\left(
A\right) \leq a_{3}.$ \ \ \ \ \ \ $\blacksquare $

Here we should observe that the only property of $M_{3}$ we used was that $%
\det M_{3}\geq 0$. Thus the conclusion of Theorem 1 is valid for some other
matrices not included in the classes $\mathcal{P}$ or $\mathcal{I}$.\vskip%
0.2in\noindent \textbf{Proof of Theorem 2.} Let $A_{r}$ be the top $r\times
r $ principal submatrix of $A$. Decompose $A_{r}$ as%
\begin{equation*}
A_{r}=D_{r}+M_{r}
\end{equation*}%
where $D_{r}$ is diagonal and $M_{r}$ off-diagonal. The matrix $iM_{r}$ is a
real skew-symmetric matrix. So, the nonzero eigenvalues of $iM_{r}$ are
purely imaginary and occur in conjugate pairs. Thus the nonzero eigenvalues
of $M_{r}$ occur in $\pm $ pairs. This shows that%
\begin{equation}
\lambda _{k}\left( M_{r}\right) \leq 0\text{ \ \ for \ }1\leq k\leq \lceil 
\frac{r}{2}\rceil .  \tag{10}
\end{equation}%
Now let $1\leq k\leq \lceil \frac{n}{2}\rceil .$ Using, successively, the
interlacing principle, Weyl's inequality and (10), we get%
\begin{equation*}
\lambda _{k}\left( A\right) \leq \lambda _{k}\left( A_{2k-1}\right) \leq
\lambda _{k}\left( M_{2k-1}\right) +a_{2k-1}\leq a_{2k-1}.
\end{equation*}%
$\blacksquare $

We now give two examples to show why for the case of matrices with
nonnegative off-diagonal entries we have to be content just with inequality
(3). Let $A$ be the $4\times 4$ matrix 
\begin{equation*}
A=\left[ 
\begin{array}{cccc}
0 & 1 & 1 & 1 \\ 
1 & 0 & 1 & 1 \\ 
1 & 1 & 0 & 1 \\ 
1 & 1 & 1 & 0%
\end{array}%
\right] .
\end{equation*}%
The $4\times 4$ matrix $E$ all whose entries are equal to one has
eigenvalues $(4,0,0,0)$. So the matrix $A=E-I$ has eigenvalues $(3,-1,-1,-1)$%
. Thus $\lambda _{3}=-1$, and the inequality (4) does not hold in this case.%
\vskip0.2in\noindent Let $B$ be the $5\times 5$ matrix%
\begin{equation*}
B=\left[ 
\begin{array}{ccccc}
0 & 0 & 1 & 1 & 0 \\ 
0 & 0 & 0 & 1 & 1 \\ 
1 & 0 & 0 & 0 & 1 \\ 
1 & 1 & 0 & 0 & 0 \\ 
0 & 1 & 1 & 0 & 0%
\end{array}%
\right] .
\end{equation*}%
Then $B=S^{2}+S^{3}$, where $S$ is the shift matrix%
\begin{equation*}
S=\left[ 
\begin{array}{ccccc}
0 & 1 & 0 & 0 & 0 \\ 
0 & 0 & 1 & 0 & 0 \\ 
0 & 0 & 0 & 1 & 0 \\ 
0 & 0 & 0 & 0 & 1 \\ 
1 & 0 & 0 & 0 & 0%
\end{array}%
\right] .
\end{equation*}%
The eigenvalues of $S$ are the fifth roots of 1. Using this one readily sees
that the eigenvalues of $B$ are $2$, $2\cos \frac{2\pi }{5}$ and $2\cos 
\frac{4\pi }{5}$, the first of these with multiplicity one and the latter
two with multiplicities two each. In particular, $\lambda _{3}>0$ and the
assertion $\lambda _{3}\leq a_{5}$ in the first inequality (5) does not hold
in this case.

\vskip0.2in\noindent \textbf{Acknowledgements}. The second author thanks
Ashoka University for arranging his visit \ during Dec 2021- Jan 2022.


\bigskip

\begin{thebibliography}{9}
\bibitem{1} A. Berman and M. Farber, \textit{A lower bound for the second
Laplacian eigenvalues of weighted graphs}, Electron. J. Linear Algebra, 
\textbf{22} (2011), 1179-1184.

\bibitem{2} R. Bhatia, \textit{Matrix Analysis}, Springer Verlag New York,
(1997).

\bibitem{3} Z. Charles, M. Farber, C. R. Johnson, L.K. Shaffer, \textit{The
relation between the diagonal entries and eigenvalues of a symmetric matrix,
based upon the sign patterns of its off-diagonal entries}, Linear Algebra
Appl; \textbf{438} (2013), 1427-1445.
\end{thebibliography}
\end{document}